\documentclass[12pt,reqno]{amsart}
\usepackage{amsfonts, amsmath, amssymb, amsthm, amsbsy}

\newtheorem{theor}{Theorem}
\newtheorem{lemma}[theor]{Lemma}
\newtheorem{cor}[theor]{Corollary}
\newtheorem{state}[theor]{Proposition}
\theoremstyle{remark}
\newtheorem{rem}{Remark}
\newtheorem{example}{Example}
\theoremstyle{definition}
\newtheorem{define}{Definition}

\newcommand{\BBC}{{\mathbb{C}}}
\newcommand{\BBR}{{\mathbb{R}}}

\newcommand{\BBE}{{\mathbb{E}}}
\newcommand{\BBZ}{{\mathbb{Z}}}

\newcommand{\const}{\mathop{\rm const}\nolimits}

\newcommand{\cI}{\mathcal{I}}
\newcommand{\cL}{\mathcal{L}}

\newcommand{\dd}{\partial}

\DeclareMathOperator{\arth}{arctanh}

\newcommand{\by}[1]{\textit{{#1}}}
\newcommand{\jour}[1]{\textit{{#1}}}
\newcommand{\vol}[1]{\textbf{{#1}}}
\newcommand{\book}[1]{\textrm{{#1}}}

\begin{document}
\rightline{ISPUmath-4/2005}

\title{On spiral minimal surfaces}

\date{February 20, 2006}

\author{A.\,V.\,Kiselev}
   %%%\thanks{A.\,K. was partially supported by University of Lecce
   %%%grant No\,650~CP/D.}

\address{Department of Higher Mathematics,\ Ivanovo State
Power University,\ Rabfakovskaya str.\,34, Ivanovo, 153003
Russia.}

\curraddr{\textup{(A.\,K.):}
Department of Physics, Middle East Technical University,
06531 Ankara, Turkey.}

\email{arthemy@newton.physics.metu.edu.tr}

\author{V.~I.~Varlamov}



\begin{abstract}
A class of spiral minimal surfaces in~$\BBE^3$
is constructed using a symmetry reduction. The new surfaces
are invariant with respect to
the composition of rotation and dilatation.
The solutions are obtained in closed form
%through the Legendre transformation
and their asymptotic behaviour is described.
\end{abstract}

\maketitle

\subsection*{Introduction}
In this paper we construct the two\/-\/dimensional minimal surfaces
$\varSigma\subset{\BBE}^3$ which are invariant with respect to
the composition of the dilatation centered at the origin and
a rotation of space around the origin.
We present the solutions in closed form through the Legendre
transformation, and we describe their asymptotic behaviour
using a non\/-\/parametric representation.
Also, we construct a special solution of the problem.

The paper is organized as follows.
In Sec.~\ref{SecSymRed} we formulate the problem of a symmetry reduction
for the minimal surface equation, and we express its general solution in
parametric form using the Legendre transformation.
In Sec.~\ref{SecRiccati} we convert the reduction problem to the
auxiliary Riccati equation and classify its solutions.
Then in Sec.~\ref{SecPhasePlane} we describe the phase portrait of
a cubic\/-\/nonlinear ODE~\cite{Bila} whose solutions determine the
profiles of the minimal surfaces on some cylinder in~$\BBE^3$.
Finally, in Sec.~\ref{SecMinSurf} we construct two classes of
the spiral minimal surfaces and indicate their asymptotic approximations.
The shape of solutions that belong to the first class
resembles the spiral galaxies. The second class of solutions is composed
by helicoidal spiral surfaces with exponentially growing
helice steps.

  %Passing from the auxiliary Riccati equation and
  %reduction~\eqref{BilaCubic} to the $\phi$-\/invariant solutions of
  %Eq.~\eqref{EMS}, we formulate the nontrivial for
  %attaching the graphs of solutions for Eq.~\eqref{BilaCubic}.
  %Next, we describe the asymptotic behaviour of the new minimal
  %surfaces~$z=\chi(\xi$,\ $\eta)$. In particular, we prove the existence
  %of the minimal surfaces that depend exponentially on the angular
  %coordinate~$\theta$ on the plane $0\xi\eta$. Also, we find a class of
  %the minimal surfaces whose shape

%The notions and concepts from geometry of differential equations
%are standard, see~\cite{ClassSym} for details.

\section{The symmetry reduction problem}\label{SecSymRed}
Let $\xi$, $\eta$, and $z$ be the Cartesian coordinates in
space~$\BBE^3$. Consider the surfaces which are locally
defined by graphs of functions,
$\varSigma=\{ z=\chi(\xi$, $\eta) \mid (\xi$,\ $\eta)\in
{\mathcal{V}}\subset{\BBR}^2\}$.
A two\/-\/dimensional surface $\varSigma$ is minimal (that is, the
minimum of area is realized by~$\varSigma$) if the function~$\chi$
satisfies the minimal surface equation
\begin{equation}\label{EMS}   %=1
(1+\chi_\eta ^2)\,\chi_{\xi\xi}-2\chi_\xi \chi_\eta \chi_{\xi\eta }
  +(1+\chi_\xi^2)\,\chi_{\eta \eta }=0.
\end{equation}
In the sequel,\label{ModuloSym}
we investigate the properties of Eq.~\eqref{EMS}
and structures related with it up to its discrete symmetry
$\xi\leftrightarrow\eta$,\ $z\mapsto-z$.

The Lie algebra of classical point symmetries for Eq.~\eqref{EMS}
is generated~\cite{Bila} by three translations %$\Phi_i$
along the respective coordinate axes,
three rotations $\Psi_j$ of the coordinate planes
around the origin, and the dilatation $\Lambda$
centered at the origin.
Let $\Psi_1$ be the (infinitesimal) rotation of the plane
$0\xi\eta$, and consider the generator $\phi=\Psi_1+\Lambda$
of a symmetry of Eq.~\eqref{EMS}.
In~\cite{Bila}, the problem of symmetry reduction for
Eq.~\eqref{EMS} by the composition $\phi$ of rotation and dilatation
was posed, although no attempts to find at least one solution were
performed. By constructions, each solution of the symmetry reduction
problem is a minimal surface invariant with respect to the
rotation of the plane $0\xi\eta$ and, simultaneously,
the dilatation of the entire space~${\BBE}^3$.

\subsection{Solutions via the Legendre transform}
The general solution of the reduction problem for
Eq.~\eqref{EMS} by its symmetry~$\phi$ is obtained~\cite{YS2005}
in parametric form using the Legendre transformation.
Now we describe a two\/-\/parametric family of the $\phi$-\/invariant
minimal surfaces.
Although, we claim that the reduction provides a special
solution that is not incorporated in this family. The present paper is
essentially devoted to the description of asymptotic properties of
the two classes of these surfaces.

First we recall that the generating section of the
symmetry at hand is
\begin{equation}\label{Symmetry}
\phi=\chi-(\xi-\eta )\chi_\xi-(\xi+\eta )\chi_\eta,
\end{equation}
see~\cite{ClassSym, Olver} for details.
The $\phi$-\/invariance condition for $\chi(\xi,\eta)$ is~$\phi=0$.
Next, consider the Legendre transformation
\begin{equation}\label{LegT}
{\cL}=\{w=\xi \chi_\xi+\eta \chi_\eta -\chi,\quad
p=\chi_\xi,\quad q=\chi_\eta \}.
\end{equation}
Now we act by the Legendre transformation~$\cL$ onto the system
composed by Eq.~\eqref{EMS} and the equation $\phi=0$.
Then from~\eqref{EMS} we obtain the linear elliptic
equation   %%%~\cite{Legendre}
\begin{subequations}\label{SystemInLeg}
\begin{gather}\label{eqGianni}  %%=3a
(1+p^2)\,w_{pp}+2pq\,w_{pq}+(1+q^2)\,w_{qq}=0,\\
\intertext{and from the invariance condition
$\phi=0$ we get the equation}
w+q\,w_p-p\,w_q=0.\label{LegSym}
\end{gather}
\end{subequations}
Let $(\varrho$, $\vartheta)$ be the polar coordinates on the
plane $0pq$ such that $p=\varrho\,\cos\vartheta$ and
$q=\varrho\,\sin\vartheta$. Then Eq.~\eqref{LegSym} acquires the form
$w-\tfrac{\partial}{\partial\vartheta} w=0$, whence it follows that
$w=\omega(\varrho)\cdot\exp(\vartheta)$.
Substituting $w(\varrho$, $\vartheta)$ in~\eqref{eqGianni},
we obtain the equation
\begin{equation}\label{RadialComponent}
\varrho^2\cdot(1+\varrho^2)\,\omega''(\varrho) +
\varrho\,\omega'(\varrho) + \omega(\varrho)=0.
\end{equation}
Its complex\/-\/valued solutions are
\begin{equation}\label{GenSolLeg}
\omega_{\pm} =
{\bigl(2+\varrho^2\bigr)}^{\frac{1}{2}}\,
\exp\Bigl(\pm\arctan{\bigl(
     %\frac{1}{\sqrt{
-(1+\varrho^2)\bigr)}^{-\frac{1}{2}}
\pm \arth{\bigl(
     %\sqrt{
-(1+\varrho^2)\bigr)}^{\frac{1}{2}}
\Bigr).
\end{equation}
Their real and imaginary parts define a basis of real solutions for
Eq.~\eqref{RadialComponent}.

The inverse Legendre transform is
${\cL}^{-1}=\{ \xi=w_p$, $\eta =w_q$, $\chi=pw_p+qw_q-w\}$.
Rewriting $\cL^{-1}$ in the polar coordinates
$\varrho$ and $\vartheta$, we finally get
\begin{equation}\label{InvLegPolar}
\xi=\Bigl[\frac{\dd\omega}{\dd\varrho}\cos\vartheta +
  \frac{\omega}{\varrho}\sin\vartheta\Bigr]\cdot\exp\vartheta,\ %
\eta=\Bigl[\frac{\dd\omega}{\dd\varrho}\sin\vartheta -
  \frac{\omega}{\varrho}\cos\vartheta\Bigr]\cdot\exp\vartheta,\ %
\chi=\Bigl[\varrho\cdot\frac{\dd\omega}{\dd\varrho}-\omega\Bigr]\cdot
   \exp\vartheta.
\end{equation}
Formulas~\eqref{InvLegPolar} provide a parametric representation
of the generic $\phi$-\/invariant minimal surfaces.
Special minimal surfaces are defined by special solutions of
Eq.~\eqref{RadialComponent} different from family~\eqref{GenSolLeg}.

\begin{rem}
Representation~\eqref{InvLegPolar} of the minimal surfaces based on
the linearizing Legendre transformation, see Eq.~\eqref{LegT},
is not a unique way to describe open arcwise connected minimal surfaces
$\varSigma=\{\vec{\xi}=\vec{\xi}(p,q)\}\in\BBE^3$ in parametric form
(e.g., we have $p=\chi_\xi$ and $q=\chi_\eta$
in formulas~\eqref{SystemInLeg}). An alternative is given through
the Ennepert\/--\/Weierstrass representation~\cite{Nitsche}
\begin{subequations}\label{Weierstrass}
\begin{gather}
\varSigma=\bigl\{\vec{\xi}(\zeta) \mid
\vec{\xi}={\vec{\xi}}_0 +
\text{Re}\,\int_0^\zeta \vec{\Phi}(\lambda)\,d\lambda,\ %
\zeta\in{\BBZ}\subset{\BBC},\ {\vec{\xi}}_0\in{\BBE}^3
\bigr\},\\
\intertext{here $\vec{\Phi}(\zeta)=\{\Phi_1(\zeta)$,\ $\Phi_2(\zeta)$,\ %
$\Phi_3(\zeta)\}$ is a triple of complex analytic functions
that satisfy the constraint}
\Phi_1^2+\Phi_2^2+\Phi_3^2=0.
\end{gather}
\end{subequations}
We see that $\text{Re}\,\zeta$ and $\text{Im}\,\zeta$ are
the new parameters on the minimal surface~$\varSigma$.
It would be of interest to obtain the Weierstrass
representation~\eqref{Weierstrass} for the $\phi$-\/invariant
minimal surfaces~\eqref{InvLegPolar}.
\end{rem}

\subsection{Reduction of Eq.~\eqref{EMS} to an ODE}
Let us return to the system composed by Eq.~\eqref{EMS} and the
constraint $\phi=0$. This symmetry reduction leads to one scalar
second order ODE with a cubic nonlinearity. Namely, we have

\begin{state}[\cite{Bila}]\label{BilaReductionState}  %=2
Let $(\rho$,\ $\theta)$ be the
polar coordinates on the plane~$0\xi\eta$
such that $(z$,\ $\rho$,\ $\theta)$ are the cylindric coordinates
in space. Consider a minimal surface $\varSigma\subset{\BBE}^3$
invariant w.r.t.\ symmetry~\eqref{Symmetry}
   %%%$\phi=\chi-(\xi-\eta )\chi_\xi-(\xi+\eta )\chi_\eta$
of Eq.~\eqref{EMS}. Then $\varSigma$ is defined by the formula
$$
z=\rho\cdot h(\theta-\ln\rho),
$$
where the function $h(q)$ of $q=\theta-\ln\rho$
satisfies the equation
\begin{equation}\label{BilaCubic}   %=2
h''\cdot(h^2+2) + h - 2h' - (h')^3 - (h-h')^3 = 0.
\end{equation}
\end{state}

\begin{cor}\label{DrawBilaSpiral}   %=1
Assume that $\varSigma$~is a minimal surface invariant
w.r.t.\ composition~$\phi$ of rotation and dilatation,
see~\eqref{Symmetry}.
Let $\Pi=\smash{\bigl\{z=h(q){\bigr|}_{\rho=1}\bigr\}}$
be the profile defined by $\varSigma$ on the
cylinder $Q=\{\rho=1\}\subset{\BBE}^3$, here~$h$ is a solution of
Eq.~\eqref{BilaCubic}.
Then the surface~$\varSigma$ is extended from~$Q$
along the intersections of the logarithmic spirals
$\rho= {\const} \cdot\exp(\theta)$, $\theta\in\BBR$,
and the cones $z=\rho\cdot h(q)$.
The $\phi$-\/invariant surface~$\varSigma$ has a singular point at
the origin.
\end{cor}

The subsitution $h=x$, $h'=y(x)$ maps Eq.~\eqref{BilaCubic} to the
cubic\/-\/nonlinear equation
\begin{equation}\label{phase_crv}   %=4
(x^2+2)y\frac{dy}{dx}={2y-x+y^3-(x-y)^3};
\end{equation}
equation~\eqref{phase_crv} defines the phase portrait of
Eq.~\eqref{BilaCubic}.

\begin{rem}\label{Italy-1917}
Consider the family of equations
\begin{equation}\label{AbelFamily}   %=5
(x^2+2)yy'=-x+\alpha y+y^3+(y-x)^3
\end{equation}
that contains~\eqref{phase_crv} as a particular case when~$\alpha=2$.
The change of variables $\smash{\frac{1}{u(x)}}=(x^2+2)y(x)$
in~\eqref{AbelFamily} leads to the Abel equation
$u'=f_3(x)u^3+f_2(x)u^2+f_1(x)u+f_0(x)$,
which is integrable in quadratures~\cite{Scalizzi} if~$\alpha=3$.
For $\alpha=2$ we have $f_0=\smash{-\frac{2}{(x^2+2)^2}}$,
$f_1=\smash{\frac{x}{x^2+2}}$, $f_2=-(3x^2+2)$, and
$f_3=x(x^2+1)(x^2+2)$.
The Abel equation can be expressed in the normal form
$\eta'=\eta^3+H(x(\xi))$, where $\eta=\eta(\xi)$, the prime
denotes the derivative w.r.t.~$\xi$, and the function $H(x)$
is constructed as follows. Let us make another change of coordinates,
$u(x)=\omega(x)\eta(\xi)-\frac{f_2}{3f_3}$, where
$\omega(x)=\smash{\exp(\int(f_1-\frac{f_2^2}{3f_3})dx)}=
\smash{\frac{ (x^2+2)^{1/6} }{ x^{2/3}(x^2+1)^{5/6} }}$ and
$\xi=\xi(x)=\int\frac{(x^2+2)^{4/3}}{x^{1/3}(x^2+1)^{2/3}}dx$.
The function~$x(\xi)$ is defined by inverting the expression
for $\xi(x)$. Then the function~$H(x)$ is obtained from the relation
$$
f_3\cdot\omega^3\cdot H(x)=f_0+{{\left(\frac{f_2}{3f_3}\right)}'_x}-
{\frac{f_1f_2}{3f_3}+\frac{2f_2^2}{27f_3^2}},
$$
whence we finally get
\begin{equation}\label{NormalTail}
H(x)=\left(\frac{x^2+1}{x^2+2}\right)^{3/2}\cdot
\frac{-20}{27x^2(x^2+1)^2(x^2+2)^2};
\end{equation}
note that
the inverse function $x(\xi)$ must be substituted for~$x$
in~\eqref{NormalTail}.
Hence we conclude that the coordinate transformation $u(x)=
\frac{ (x^2+2)^{1/6}}{ x^{2/3}(x^2+1)^{5/6} }\cdot\eta(\xi)+
\frac{3x^2+2}{3x(x^2+2)(x^2+1)}$
maps the Abel equation to its normal form $\eta'=\eta^3+H(x(\xi))$.
\end{rem}

\subsection{Exact solutions of the limit analogue
   of Eq.~\eqref{phase_crv}}\label{SecLimitEq}
In this subsection we consider the limit analogue of
system~\eqref{phase_crv} whose right\/-\/hand side contains
the terms that provide maximal contribution as
$R=\sqrt{x^2+y^2}\to\infty$ on~$0xy$.
To this end, we note that for large~$R$ the right\/-\/hand side of
Eq.~\eqref{phase_crv}
 %$$ \tag{ } %\tag{\ref{phase_crv}}
 %\frac{dy}{dx}=\frac{2y-x+y^3-(x-y)^3}{y(x^2+2)}
 %$$ %\end{equation}
is equivalent to its ``cubic'' component
$\bigl[{y^3-(x-y)^3}\bigr]\cdot{y^{-1}\cdot x^{-2}}$.
Let us find exact solutions of the equation
\begin{equation}\label{phase_crv3}   %=26
\frac{dy}{dx}=\frac{y^3-(x-y)^3}{yx^2}.
\end{equation}
   %as~$R\to\infty$.
   %Equation~\eqref{phase_crv3} can be explicitly integrated.
Using the substitution~$s=\frac{y}{x}$ in~\eqref{phase_crv3},
we obtain the equation
$$
x\,\frac{ds}{dx}=\frac{2s^3+3s-4s^2-1}{s}.
$$
Note that $2s^3-4s^2+3s-1=(s-1)(2s^2-2s+1)$.
Therefore we get
\begin{equation}\label{ApproxDelta}
\int\frac{s\,ds}{(s-1)(2s^2-2s+1)}=\frac12\ln\frac{(s-1)^2}{2s^2-2s+1}
=\ln|\delta x|,
  %\int\frac{s\,ds}{(s-1)(2s^2-2s+1)}=
  %\int\frac{ds}{s-1}-\int\frac{2s-1}{2s^2-2s+1}ds=
\end{equation}
whence we finally obtain the integral
$\frac{(s-1)^2}{2s^2-2s+1}=\delta x^2$,
here~$\delta\geq 0$.

Suppose $\delta=0$, then we have~$s=1$. Thus we get the solution
$y=x$ of approximation~\eqref{phase_crv3}.
In Proposition~\ref{AsymptSepState} below we prove that the diagonal
is the asymptote for a solution of Eq.~\eqref{phase_crv}.
  %%% simultaneously, this is the asymptote for the separatrix~$y_2$.
If $\delta>0$, then we solve the quadratic equation w.r.t.\ $s$ and
obtain
$$
s=\frac{1-\delta x^2\pm\sqrt{\delta }\cdot|x|\cdot\sqrt{1-\delta x^2}}%
   {1-2\delta x^2}.
$$     %%% s=\frac{y}{x};
Hence we conclude that
\begin{subequations}\label{twosol}   %=27
\begin{align}
\bar y_1
&=\frac{x\cdot\sqrt{1-\delta x^2}}{1-2\delta x^2}(\sqrt{1-\delta x^2}
    +\sqrt{\delta }\cdot|x|),   \label{sol1}\\ %=27a
\bar y_2
&=\frac{x\cdot\sqrt{1-\delta x^2}}%
   {\sqrt{1-\delta x^2}+\sqrt{\delta }\cdot|x|}.  \label{sol2} %=27b
\end{align}
\end{subequations}
Solutions~\eqref{twosol} are defined for $|x|<\frac1{\sqrt{\delta}}$,
$|x|\ne\frac1{\sqrt{2\delta}}$.
Suppose $x>0$, then we conclude that
$$
\lim_{x\to\frac1{\sqrt{2\delta }}-0}\bar y_1(x)=+\infty.
$$
Finally, consider solution~\eqref{sol1} with the Cauchy data
$(x_0,y_0)$, $y_0\ne x_0$, located on a large circle.
Then the above reasonings yield that
$$
\delta =\frac{|\frac{y_0}{x_0}-1|}{\sqrt{2y_0^2-2y_0x_0+x_0^2}}.
$$

We claim that the behaviour of solutions of
initial equation~\eqref{phase_crv} is correlated with the solutions of
limit equation~\eqref{phase_crv3}, see Remark~\ref{CompareWithLimit}
on p.~\pageref{CompareWithLimit}.

%%%%%%%%%%%%%%%%%%%%%%%%%%%%%%%%%%%%%%%%%%%%%%%%%%%%%%%%%%%%%%%%%%%%%%%%
\section{Solutions of the auxiliary Riccati equation}\label{SecRiccati}
In this section, we start to investigate the phase portrait of
Eq.~\eqref{BilaCubic}. We show that the phase curves are
described by solutions of the auxiliary Riccati
equation~\eqref{Riccati}. Using the topological method
by Warzewski~\cite{Hartman}, we conclude that Eq.~\eqref{Riccati} has
one unstable solution and a large class of stable solutions.
On the phase plane for Eq.~\eqref{BilaCubic} these solutions are
represented by a pair of centrally symmetric separatrices and by
trajectories that repulse from the separatrices and have vertical
asymptotes at infinity, respectively. This will be discussed in
the next section.

\label{SecTransform}
Consider the autonomous system associated with~\eqref{phase_crv},
\begin{equation}\label{xysyst}   %=6
\frac{dx}{dt}=2y+yx^2,\quad \frac{dy}{dt}=2y-x+y^3-(x-y)^3.
\end{equation}
Integral curves for Eq.~\eqref{phase_crv} are assigned
to integral trajectories of system~\eqref{xysyst}.
Consider the linear change of variables $x=2(z_1-z_2)$, $y=2z_1$.
Then system~\eqref{xysyst} is transformed to
\begin{equation}\label{z_syst}   %=7
\frac{dz_1}{dt}=z_1+z_2+4(z_1^3+z_2^3),\quad
\frac{dz_2}{dt}=-z_1+z_2+4z_2(2z_1^2-z_1z_2+z_2^2).
\end{equation}
The linear approximation for~\eqref{z_syst} has an unstable focus at
the origin; indeed, the eigenvalues are $\lambda_{1,2}=1\pm i$.
By a Lyapunov's theorem, nonlinear systems~\eqref{xysyst}
and~\eqref{z_syst} will exhibit analogous behaviour near the origin.

Next, we use the substitution $z_1=r\cos\varphi$, $z_2=r\sin\varphi$
in~\eqref{z_syst},
 %\footnote{Замечание.
 %Народ единогласно восстал против матрицы линейной подкрутки,
 %которую следует поделить на 2, и радиального уравнения, из которого
 %тем же методом изгоняется коэффициент 4 -- чтобы $\dot{r}=r+r^3$.}
whence after some transformations we obtain the
\emph{triangular} system
\begin{subequations}\label{rfsyst}\label{polar_syst}   %=8
\begin{align}
\dot r&=r+4r^3,\label{radial}\\   %=8a
\dot{\varphi}&=-1+4r^2\sin\varphi\,(\cos\varphi-\sin\varphi).
   \label{angular}   %=8b
\end{align}
\end{subequations}
For any $r_0>0$ the Cauchy problem for Eq.~\eqref{radial} has the
solution
\begin{equation}\label{rsol}  %=9
r(t,r_0)=\frac{r_0\exp(t)}{\sqrt{1+4r_0^2\bigl(1-\exp(2t)\bigr)}}
\end{equation}
whenever~$r(0,r_0)=r_0$.
Assume that $r=r_0>0$ at $t=0$ for a solution
$\{r(t)$,\ $\varphi(t)\}$.
From~\eqref{rsol} it follows that any solution
of system~\eqref{polar_syst} achieves the infinity
on the plane~$0xy$ at the finite time
$t^*=\smash{\frac12\ln(1+\frac1{4r_0^2})}$, and we see that
$t^*\to\infty$ if $r_0\to+0$.
Obviously, we have $\dot\varphi\to-1$ for $r\to+0$,
that is, the trajectories in a small neighbourhood of the origin
are the spirals that unroll clockwise.
These reasonings describe the behaviour of the trajectories of
systems~\eqref{xysyst} and~\eqref{z_syst}.
   %%%Now we formulate the assertion.

\begin{state}\label{Transversality}  %=3
The trajectories ${z_1}(t)$,\ $z_2(t)$ are transversal to any circle
centered at the origin. Therefore the phase curves $\{x(t)$, $y(t)\}$
for Eq.~\eqref{BilaCubic} are transversal to any cetrally symmetric
ellipsis %in $0xy$
that corresponds to a circle on the plane~$0z_1z_2$
\textup{(}the major axes of these ellipses are located along the
diagonal~$y=x$\textup{)}. Hence we deduce that system~\eqref{xysyst}
has no cycles and equation~\eqref{BilaCubic} has no periodic solutions
except zero, which corresponds to the stable point on the
phase plane.
A unique stable point $(0,0)$ of system~\eqref{xysyst} is an
unstable focus. The spiral\/-\/type phase curves for
Eq.~\eqref{BilaCubic} unroll clockwise around this point.
Any trajectory of system~\eqref{xysyst} that does not coincide with the
origin achieves the infinity at a finite time.
\end{state}

\begin{rem}
The equations in system~\eqref{xysyst} do not change under the
involution $(x,y)$ $\leftrightarrow(-x,-y)$. Hence for any
trajectory of systems~(\ref{xysyst}--\ref{polar_syst}) its
image under the central symmetry is also a trajectory.
Further on, we investigate properties of the phase curves
$\{x(t)$,\ $y(t)\}$ up to this symmetry
(or, equivalently, up to the transformation $\varphi\mapsto\varphi+\pi$
of the angular coordinate~$\varphi$).
\end{rem}

\subsection{} %Подстановка тангенс. Риккати и Шрёдингер.
In this subsection we study the behavoiur of solutions of
system~\eqref{polar_syst} as~$r\to+\infty$.
Let us divide Eq.~\eqref{angular} by Eq.~\eqref{radial}.
Thus we obtain %the equation
\begin{equation}\label{dfdr}   %=10
\frac{d\varphi}{dr}=
   \frac{-1+4r^2\sin\varphi(\cos\varphi-\sin\varphi)}{r+4r^3}.
\end{equation}
Substituting $\tau=4r^2$ in~\eqref{dfdr}, we get the equation
\begin{equation}\label{dfdt}   %=11
\frac{d\varphi}{d\tau}=
   \frac{-1+\tau\sin\varphi(\cos\varphi-\sin\varphi)}{2\tau(1+\tau)}.
\end{equation}
After some trigonometric simplifications in the r.h.s.\ of
Eq.~\eqref{dfdt} we arrive at
\begin{equation}\label{dfdtau}\tag{\ref{dfdt}${}'$}
\frac{d\varphi}{d\tau}=\frac{(\sqrt2-1)\tau-2
   -2\sqrt2\tau\sin^2(\varphi-\frac{\pi}8)}
{4\tau(\tau+1)}.
\end{equation}
Next, divide both sides in Eq.~\eqref{dfdtau} by
$\cos^2(\varphi-\frac{\pi}8)$ assuming that
$\varphi\ne-\frac{3\pi}8+\pi n$, $n\in{\BBZ}$.
Also, we put $u=\tan(\varphi-\frac{\pi}8)$ by definition.
Hence we finally obtain the Riccati equation
\begin{subequations}\label{Riccati}   %=13
\begin{gather}
\frac{du}{d\tau}=b(\tau)-a(\tau)u^2,\\
\intertext{where}
a(\tau)=\frac{(\sqrt2+1)\tau+2}{4\tau(\tau+1)},\quad
b(\tau)=\frac{(\sqrt2-1)\tau-2}{4\tau(\tau+1)}.
\end{gather}
\end{subequations}
The Riccati equation~\eqref{Riccati} can be further
transformed~\cite{Zelinkin} to the
Schr\"odinger equation with the potential~$b/a$ and zero energy.

\subsection{The Warzewski theorem}
First let us recall some definitions~\cite{Hartman}.

\begin{define}
Let $\Omega^0\subset\Omega$ be an open subset of a domain~$\Omega$.
Consider the Cauchy problem
\begin{subequations}\label{Cauchy}   %=15
\begin{align}
y'=f(t,y), \label{syst}\\   %=15a
y(t_0)=y_0.    \label{nu}   % \tag{15b}
\end{align}
\end{subequations}
Suppose $y=y(t)$ is a solution of~\eqref{Cauchy}.
A point $(t_0,y_0)\in\Omega\cap\partial\Omega^0$ is called
an \emph{exit point} with respect to equation~\eqref{syst}
and the domain~$\Omega^0$ if for any solution~$y(t)$
satisfying~\eqref{nu} there is a constant $\varepsilon>0$ such that
$(t,y(t))\in\Omega^0$ for all~$t$ in the interval
$t_0-\varepsilon\leq t<t_0$.
An exit point $(t_0,y_0)$ for the domain $\Omega^0$ is
a \emph{strict exit point} %from~$\Omega^0$
if $(t,y(t))\notin\overline{\Omega^0}$ whenever
$t_0<t\leq t_0+\varepsilon$ for some $\varepsilon>0$.
We denote by~$\Omega^0_e$ the set of all exit points for the
domain~$\Omega^0$, and let
$\Omega^0_{se}$ denote the set of all strict exit points.
\end{define}

Now we analyze the behaviour of solutions of Eq.~\eqref{Riccati}
as $\tau\to+\infty$. First we note that
$\tan\frac{\pi}8=\sqrt2-1$; hence we put $u_1=\sqrt2-1$,
$u_2=1-\sqrt2$. The above notation corresponds to the angles
$\varphi_1=\frac{\pi}4$ and $\varphi_2=0$, respectively.
Also, we set $\tau_0=2(\sqrt2+1)+\delta_0$, where~$\delta_0>0$
is arbitrary.

Further, we introduce two closed domains $D_1$ and $D_2$ in the right
half\/-\/plane $\tau>0$ of the coordinate plane~$(\tau;u)$: we let
\begin{equation}\label{domain1} %=14 %\label{domain2}
\begin{aligned}
D_1&=\{(\tau;u)|\tau\geq\tau_0,\ 0\leq u\leq u_1\},\\
D_2&=\{(\tau;u)|\tau\geq\tau_0,\ u_2\leq u\leq0\}.
\end{aligned}
\end{equation}
Next, we note that the equality $\frac{du}{d\tau}=0$ is valid on the
curves $u(\tau)=\pm\sqrt{{b(\tau)}\bigr/{a(\tau)}}$. Therefore we
define the third domain $D_3\subset D_1\cup D_2$ through
$$
D_3=\left\{(\tau,u)|\tau>2(\sqrt{2}+1),\ %
-\sqrt{{b(\tau)}\bigr/{a(\tau)}}
<u<
\sqrt{{b(\tau)}\bigr/{a(\tau)}}\right\}.
$$
It is easy to check that the inequality $\frac{du}{d\tau}>0$
holds in the domain~$D_3$.

By definition, put $f(u,\tau)=b(\tau)-a(\tau)u^2$.
Let us describe the inclination field for Eq.~\eqref{Riccati}
on the lines $u=u_1$, $u=0$, and $u=u_2$. We have
$f(u_1,\tau)=f(u_2,\tau)=-\frac{1+(\sqrt2-1)^2}{2\tau(\tau+1)}<0$
whenever $\tau\geq\tau_0$, and we also have
$f(0,\tau)=b(\tau)>0$ under the same assumption.
This argument shows that
all points on the lines $u=0$ and $u=u_1$ are the strict entry
points with respect to~$D_1$. Simultaneously, the lines
$u=0$ and $u=u_2$ are composed by the strict exit points
w.r.t.\ the closed domain~$D_2$.

The Warzewski theorem and Example~\ref{ExampleWarz} below are borrowed
from~\cite{Hartman}. Using them, we prove the existence of solutions of
Eq.~\eqref{Riccati} that do not leave the domains $D_1$ and~$D_2$.

\begin{theor}[\cite{Hartman}]\label{WarzTh}   %=2  с.~332--334
Let $f(t,y)$ be a continuous function on an open set $\Omega$
of points $(t$,\ $y)$, and assume that solutions of
system~\eqref{syst}
are uniquely defined by the initial condition~\eqref{nu}.
Also, let $\Omega^0\subset\Omega$ be an open subset such that
$\Omega^0_e=\Omega^0_{se}$.
Further let~$S\subset\Omega^0$ be a nonempty subset such that
\begin{itemize}
\item
the intersection $S\cap\Omega^0_e$ is a retract of~$\Omega^0_e$, but
\item
the intersection $S\cap\Omega^0_e$ is not a retract of~$S$.
\end{itemize}
Then there is at least one point $(t_0,y_0)\in S\cap\Omega^0$
such that the graph of solution~$y(t)$ of the Cauchy
problem~\eqref{Cauchy} is contained in $\Omega^0$ on its maximal right
interval of definition.
\end{theor}

\begin{example}[\cite{Hartman}]\label{ExampleWarz}
Suppose that $y$ is real and the function~$f(t,y)$ in
system~\eqref{syst} is continuous on the set $\Omega$ which coincides
with the whole plain $(t,y)$. Let $\Omega^0$ be the strip
$|y|<b,\ -\infty<t<\infty$.
The boundary of~$\Omega^0$ is contained in~$\Omega$ and
consists of the two lines~$y=\pm b$. Assume $f(t,b)>0$ and $f(t,-b)<0$
such that $\Omega^0_e=\Omega^0_{se}=\partial\Omega^0\cap\Omega$.
Next, let $S$ be the segment
$S=\{(t,y)\mid t=0,\ |y|\leq b\}$.
Then $S\cap\Omega^0_e$ consists of the two points
$(0,\pm b)$; the intersection is a retract of of the set~$\Omega^0_e$
but is not a retract of~$S$.
Theorem~\ref{WarzTh} yields the existence of a point
$(0,y_0)$, $|y_0|<b$ such that there is the solution of
the Cauchy problem for
system~\eqref{syst} with the initial condition $y(0)=y_0$.
This solution satisfies the inequality $|y(t)|<b$ for all~$t\geq0$.
\end{example}

From Example~\ref{ExampleWarz} that illustrates
Theorem~\ref{WarzTh} we obtain %the following conclusion.

\begin{cor}\label{FromWarz} %=3
There is a solution $u=\psi_1(\tau)$ of the Riccati
equation~\eqref{Riccati} that does not leave the domain~$D_1$,
which is defined in~\eqref{domain1}, for all~$\tau\geq\tau_0$.
Analogously, there is a solution $u=\psi_2(\tau)$ in $D_2$ that does
not leave~$D_2$.
\end{cor}

\begin{rem}
The Warzewski theorem guarantees the existence of a solution
$\psi_1$ in $D_1$. We claim that there are infinitely many solutions of
this class in our case. The asymptotic expansions of these solutions as
$r\to\infty$ are specified in~\eqref{psi1}. The convergence of the
expansions can be rigorously proved,
see Remark~\ref{RemMajorantSeries} on p.~\pageref{RemMajorantSeries}.
Also, we claim that the solution~$\psi_2$ in the domain~$D_2$ is unique
and unstable.
\end{rem}

Now we calculate the limits of solutions of
the Riccati equation~\eqref{Riccati} as $\tau\to\infty$.

\begin{lemma}\label{CommonLimitAboveLemma}   %=1
Let $u(\tau)$ be a solution of Eq.~\eqref{Riccati}
which is greater than~$u=\psi_1(\tau)$ for all $\tau\geq\tau_0$.
Then we have
\begin{equation}\label{LimAll}   %=16
\lim_{\tau\to\infty}(u(\tau)-\psi_1(\tau))=0.
\end{equation}
\end{lemma}

\begin{proof}
Assume the converse,
$\inf\limits_{\tau\geq\tau_0} (u(\tau)-\psi_1(\tau))=\Delta_0$,
where $\Delta_0>0$.
Then the difference $w(\tau)=u(\tau)-\psi_1(\tau)$
of two solutions for Eq.~\eqref{Riccati} satisfies the equation
\begin{equation}\label{fordiffw}   %=17
\frac{dw}{d\tau}=-a(\tau)\bigl(u(\tau)+\psi_1(\tau)\bigr)\cdot w.
\end{equation}
Integrating Eq.~\eqref{fordiffw}, we obtain
\begin{equation}\label{diffw} %=18
w(\tau,\tau_0)=w(\tau_0)\cdot\exp(-\int_{\tau_0}^{\tau}a(s)(u(s)+\psi_1(s))ds.
\end{equation}
Further recall that
$$
u(s)+\psi_1(s)\geq u(s)-\psi_1(s)\geq\Delta_0
$$
whenever $s\geq\tau_0$ and $a(s)>0$. Therefore,
\begin{multline*}
-\int_{\tau_0}^{\tau}a(s)(u(s)+\psi_1(s))ds\leq
  -\Delta_0\int_{\tau_0}^{\tau}a(s)ds=
-\frac{\Delta_0}4
\left((\sqrt2-1)\ln\frac{\tau+1}{\tau_0+1}+2\ln\frac{\tau}{\tau_0}\right).
\end{multline*}
Consequently, we have
$\lim\limits_{\tau\to+\infty}w(\tau,\tau_0)=0$.
This conclusion contradicts the assumption~${\Delta_0>0}$.
\end{proof}

Similarly we prove that a solution $u(\tau)$ tends to
$u=\psi_1(\tau)$ as $\tau\to+\infty$ if $u(\tau)$ enters the domain
$D_1$ and there it is located under the graph $u=\psi_1(\tau)$.

\begin{lemma}
The following equality holds\textup{:}
$$\lim_{\tau\to+\infty}\psi_1(\tau)=
\lim_{\tau\to+\infty}\sqrt{{b(\tau)}\bigr/{a(\tau)}}=\sqrt2-1=u_1.$$
\end{lemma}

\begin{proof}
Suppose $\psi_1(\tau)\subset D_1\cap D_3$ for all~$\tau>\tau_0$.
Then the solution $\psi_1(\tau)$ grows and is bounded. Therefore
there is the limit
$$\lim_{\tau\to+\infty}\psi_1(\tau)=d_0.$$
Next, let $u=u(\tau)$ be a solution that enters the domain
$D_3\cap D_1$ through the line $u=u_1$ at a large~$\tau_1$ and remains
in it for all $\tau>\tau_1$.
If $d_0<\sqrt2-1$, then the limit of the solution $u(\tau)$ is
$d_1>d_0$, but the limits $d_0$ and $d_1$ must coincide
by~\eqref{LimAll}. Consequently, $d_0=d_1=\sqrt2-1$.
\end{proof}

\begin{lemma}
The solution $u=\psi_2(\tau)$ remains in the domain $D_2\setminus D_3$
for all $\tau\geq\tau_0$, and its limit at infinity is
$$\lim_{\tau\to+\infty}\psi_2(\tau)=1-\sqrt2.$$
\end{lemma}

\begin{proof}
If the solution $\psi_2(\tau)$ enters the domain $D_2\cap D_3$ at some
$\tau_1>\tau_0$ and remains there, then it grows and is bounded.
Consequently, there is the limit
$\lim\limits_{\tau\to+\infty}\psi_2(\tau)=d_2$, where
$d_2\in(1-\sqrt2;0]$.
Hence for a large $\tau_2$ and $\tau\geq\tau_2$ we obtain
\[   %%%\begin{equation}\label{psisum}  %=19
\psi_1(\tau)+\psi_2(\tau)\geq\frac{\sqrt2-1+d_2}2.
\]   %%%\end{equation}
By definition, put $w(\tau)=\psi_1(\tau)-\psi_2(\tau)$.
Then from Eq.~\eqref{Riccati} it follows that
\[   %%%\begin{equation}\label{dwdtau}   %=20
\frac{dw}{d\tau}=-a(\tau)(\psi_1(\tau)+\psi_2(\tau))w,
\]   %%%\end{equation}
whence we deduce
\begin{multline}
w(\tau)=w(\tau_2)\cdot\exp\left(-\int_{\tau_2}^{\tau}a(s)(\psi_1(s)+
\psi_2(s))ds\right)
\leq\\
\leq w(\tau_2)\exp\left(-\frac{\sqrt2-1+d_2}2\cdot
\int_{\tau_2}^{\tau}a(s)ds\right).\label{psi2NotGrows}
\end{multline}
Yet we see that the r.h.s.\ in~\eqref{psi2NotGrows}
tends to $0$ as $\tau\to+\infty$. Therefore
the l.h.s.\ in~\eqref{psi2NotGrows}
must also tend to zero which is impossible.
Hence the graph $\psi_2(\tau)$ is contained in $D_2\setminus D_3$ for
large~$\tau$, that is, the solution $\psi_2$ decreases and is bounded
from below. This argument shows that
$$\lim_{\tau\to+\infty}\psi_2(\tau)=
\lim_{\tau\to+\infty}\left(-\sqrt{{b(\tau)}\bigr/{a(\tau)}}\right)
   =-(\sqrt2-1)=1-\sqrt2.
$$
This completes the proof.
\end{proof}

Analogously, suppose the graph of a solution enters the domain~$D_3$
and is located above the graph $u=\psi_2(\tau)$. Then it tends to the
solution $u=\psi_1(\tau)$. The proof is straightforward.
A solution which is less than ${u=\psi_2(\tau)}$ achieves $-\infty$
at a finite time.
Using Eqs.~\eqref{fordiffw} and~\eqref{diffw}, it can be proved that
the solution $u=\psi_2(\tau)$ in the domain $D_2\setminus D_3$ is
unique. The proof is by \emph{reduction ad absurdum}.

\subsection{Asymptotic expansions of
   the solutions $\psi(\tau)$}\label{SecInTau}
In this subsection we use the method of undetermined coefficients
and obtain the asymptotic expansion in $\tfrac{1}{\tau}$
for the solution $\psi_1^*(\tau)$ of Eq.~\eqref{Riccati}
that tends from above to $u_1=\tan\frac{\pi}{8}$ as $\tau=4r^2\to+\infty$.
Also, we get the expansion for $\psi_2(\tau)$
that tends to $u_2=-\tan\frac{\pi}{8}$ at infinity.
  %%Recall that $u=\psi^*_1(\tau)$ is the solution
  %%that tends to $u_1=\sqrt{2}-1$ as ,
  %%and the solution $u=\psi_2(\tau)$ tends to $u_2=1-\sqrt2$ as
  %%$\tau\to+\infty$.
  %%we obtain the expansions in negative powers of~$\tau$ for
  %%the solutions $\psi^*_1(\tau)$ and~$\psi_2(\tau)$.
  %%see Proposition~\ref{PsiAsympt} below.
We emphasize that the expansion for
$\psi_1^*$ provides a solution different from $\psi_1(\tau)$, which
exists by Corollary~\ref{FromWarz}. Indeed, the new function~$\psi_1^*$
tends to the limit $u_1=\sqrt2-1$ monotonously \emph{descreasing}.

\begin{state}\label{PsiAsympt}
The solutions $\psi^*_1(\tau)$ and $\psi_2(\tau)$ admit the following
asymptotic expansions\textup{:}
\begin{subequations}\label{psir2} %=21
\begin{align}
\psi^*_1(\tau)&
=\sqrt2-1+\frac{2\sqrt2(\sqrt2-1)}{\tau}+O\bigl(\frac1{\tau^2}\bigr),
    \label{psi1star}   %=21a
  \\
\psi_2(\tau)&=1-\sqrt2+\frac{2\sqrt2(\sqrt2-1)}{3\tau}
   +O\bigl(\frac1{\tau^2}\bigr). \label{psi2}    %=21b
\end{align}
\end{subequations}
The inequalities $\psi^*_1(\tau)>\sqrt2-1$ and
$\psi_2(\tau)>1-\sqrt2$ hold for large~$\tau$.
\end{state}

\begin{rem}\label{RemMajorantSeries}
Using the method of majorant series~\cite{VarlamovConverg,
VarlamovConvergISPU}, we prove that bounded solutions of the Riccati
equation~\eqref{Riccati} are assigned to expansions~\eqref{psir2}.
These solutions are real analytic in a neighbourhood of the infinity,
and the radius of this neighbourhood can be estimated.
\end{rem}

\subsection{The general case: $u=\psi(r)$.}\label{GeneralCase}
Now we construct the asymptotic expansions for
\emph{all} solutions of~\eqref{Riccati}
that tend to $\pm(\sqrt2-1)$.
This time we use the expansions in~$r$ but not in~$\tau=4r^2$.

Let us re\/-\/write the Riccati equation~\eqref{Riccati} in the form
\begin{equation}\label{one}   %=22
\frac{du}{dr}=
  \frac{2(\sqrt2-1)r^2-1}{r(4r^2+1)}
  -\frac{2(\sqrt2+1)r^2+1}{r(4r^2+1)}\cdot u^2.
\end{equation}
Its right\/-\/hand side is real analytic in $\frac{1}{r}$
if~$r>\frac{1}{2}$. Suppose $r\geq1$.
Let us use the expansion
\begin{equation}\label{two}   %=23
u(r)=w_0+\frac{w_1}{r}+\frac{w_2}{r^2}+O\bigl(\frac1{r^3}\bigr).
\end{equation}
Substituting~\eqref{two} for~$u$ in~\eqref{one}, we get
\begin{multline*}
-\frac{w_1}{r^2}-\frac{2w_2}{r^3}-\ldots=
\frac{\sqrt2-1}{2r}\cdot
 \bigl(1-\frac1{4r^2}+\frac1{16r^4}-\ldots\bigr)-{}\\
{}-\frac1{4r^3}\cdot\bigl(1-\frac1{4r^2}+\frac1{16r^4}-\ldots\bigr)-
\bigl(w_0^2+\frac{2w_0w_1}{r}+\frac{w_1^2+2w_0w_2}{r^2}+
   \ldots\bigr)\cdot{}\\
{}\cdot\Bigl(\frac{\sqrt2+1}{2r}
\cdot\bigl(1-\frac1{4r^2}+\frac1{16r^4}-\ldots\bigr)+
\frac1{4r^3}\cdot\bigl(1-\frac1{4r^2}+\frac1{16r^4}-\ldots\bigr)\Bigr).
\end{multline*}
Equating the coefficients of $1/r$ and $1/r^2$, we obtain
$$
\sqrt2-1-(\sqrt2+1)\cdot w_0^2=0,\quad w_1=(\sqrt2+1)\cdot w_0w_1.
$$
The former equation has two roots,
$w_{0,1}=\sqrt2-1$ and $w_{0,2}=1-\sqrt2$.
First we let $w_0=w_{0,1}$; then the root of the second equation is an
arbitrary real number! In this case, we use the notation~$w_1=C$.
Secondly, suppose $w_0=w_{0,2}$; then a unique root of the second
equation is~$w_1=0$.
Now we equate the coefficients of $\frac1{r^3}$, whence we get
$$
-2w_2=
-\frac{\sqrt2-1}8-\frac14-\frac{w_0^2}4
-\frac{\sqrt2+1}2\cdot(w_1+2w_0w_2)+
\frac{(\sqrt2+1)w_0^2}8.
$$
If $w_0=\sqrt2-1$, then
$w_2=\frac{\sqrt2(\sqrt2-1)}2+\frac{\sqrt2+1}2\cdot C^2$.
Alternatively, if $w_0=1-\sqrt2$, then
$w_2=\frac{\sqrt2(\sqrt2-1)}6$.
This implies that for $C=0$ we obtain the two asymptotic expansions for
the solutions $\psi_1^*(\tau)$ and $\psi_2(\tau)$, which depend on
even powers of~$r$ and which were previously found in~\eqref{psir2}.
Now we see that all other expansions involve~$r=\frac{\sqrt{\tau}}2$
explicitly.

By definition, put $D=\tfrac{1}{2}\bigl[
\sqrt{2}(\sqrt{2}-1)+(\sqrt{2}+1)\cdot C^2\bigr]$.
Then for any~$u(r)$ such that $\lim_{r\to+\infty}u(r)=u_1$
we finally have
\begin{equation}\label{psi1}   %=24
u(r,C)=\sqrt2-1+\frac{C}{r} + \frac{D}{r^2}+O\bigl(\frac1{r^3}\bigr),
\qquad C\in\BBR.
\end{equation}
An analogue of Remark~\ref{RemMajorantSeries} is also true
for~\eqref{psi1}: using the method of majorant
series~\cite{VarlamovConverg}, one readily proves the convergence of
expansion~\eqref{psi1} for large~$r$
to analytic solutions of equation~\eqref{one}, and it is also possible
to estimate the radius of convergence.

%\subsection{}
Finally, we formulate the assertion about the behaviour of solutions
of the Riccati equation~\eqref{Riccati}.

\begin{theor}
If $\tau\in(0,2[\sqrt{2}+1])$, then all solutions of Eq.~\eqref{Riccati}
decrease monotonously.
Suppose $\tau\geq\tau_0$. Then there is the unstable solution
$\psi_2(\tau)$ that tends from above to the limit
$u_2=-\tan\frac{\pi}{8}$ as $\tau\to\infty$\textup{;}
the asymptotic expansion for~$\psi_2$ is given in~\eqref{psi2}.

In the domain $D_1$, see~\eqref{domain1}, there are infinitely
many growing solutions $\psi_1(r)$ that tend to
$u_1=\tan\frac{\pi}{8}$ as $r\to\infty$.
These solutions correspond to $C<0$ in expansion~\eqref{psi1}.

All solutions which are located between~$\psi_2(\tau)$
and~$\psi_1(r)$ tend to $u_1$ as $\tau=4r^2\to\infty$.
All solutions which are located under~$\psi_2$ repulse
from it\textup{;} they decrease and achieve $-\infty$ at a finite time.

The solution $\psi_1^*(\tau)$, which is located above the
domain~$D_1$, decreases and tends to $u_1$ as
$\tau\to\infty$\textup{;} its expansion is described
by~\eqref{psi1star} \textup(or by formula~\eqref{psi1} with
$C=0$\textup{)}.
There are infinitely many other solutions $\psi_1(r)$
which are greater than~$\psi_1^*(\tau)$. These solutions
also tend to $u_1$ as $\tau\to\infty$, and their expansions at
infinity are given through~\eqref{psi1} with~$C>0$.
\end{theor}

\section{The phase plane~$0xy$}\label{SecPhasePlane}
In this section we describe the behaviour of phase curves for
Eq.~\eqref{phase_crv}. Inverting the transformations introduced on
p.~\pageref{SecTransform} and preserving the subscripts of respective
functions, we pass from solutions $\psi(\tau)$ of Eq.~\eqref{Riccati}
to solutions $\{x(t)$,\ $y(t)\}$ of system~\eqref{xysyst}, and next
we obtain solutions $y(x)$ of Eq.~\eqref{phase_crv}.
For example, the solution $\psi_1^*(\tau)$ is transformed
to~$y_1^*(x)$, and the separatrix~$y_2$ is assigned to the unstable
solution~$\psi_2$.

The exposition goes along the time~$t$ in~\eqref{xysyst}:
first we consider a neighbourhood of the origin, next we analyze the
behaviour of the trajectories at a finite distance from the origin, and
finally we describe their asymptotic expansions at the infinity
($x^2+y^2\to\infty$).

\subsection{}
From Eq.~\eqref{radial} it follows that system~\eqref{xysyst} has an
unstable focus at the origin. In a small neighbourhood of this point,
all trajectories unroll clockwise.
By Proposition~\ref{Transversality}, each trajectory is transversal to
the ellipses centered at the origin (we recall that these ellipses
correspond to the circles $r=\const$ on the plane~$0z_1z_2$);
the trajectories achieve the infinity at a finite time.

\begin{state}\label{ExtremaState}
All extrema of the phase curves for system~\eqref{xysyst}
are located on the straight line~$y=x/2$.
\end{state}

\begin{proof}
Solving the equation $P_3(x,y)=(y-x)^3+y^3+2y-x=0$ with respect to~$y$,
we obtain the real root~$x/2$.
The quotient $P_3(x,y)/(y-x/2)=2y^2-2xy+2x^2+2$ has no real roots
for~$y$ at any~$x$.
\end{proof}

Consider the phase trajectories $\{x(t)$,\ $y(t)\}$
of Eq.~\eqref{phase_crv} in a neighbourhood of the axis~$0x$.
The parts of the trajectories near the points $(x,0)$
describe the graphs of solutions~$h(q)$
of Eq.~\eqref{BilaCubic} near their extrema.
The trajectories are approximated with the circles
of radius $\varrho(x)=|x|\cdot\frac{x^2+1}{x^2+2}$ centered
at the points $(x_0$,\ $0)$, where $x_0=\frac{|x|}{x^2+2}$.
Obviously, we have $x_0\to0$ and $\varrho(x)\to|x|$ as~$|x|\to\infty$.

\begin{rem}[On inflection points]
In the first and second quadrants of~$0xy$ (and owing to the central
symmetry of the phase portrait, in the third and fourth quadrants,
respectively) there is the curve~$\cI$ that consists of the inflection
points of the phase curves. The curve~$\cI$ is composed by two
components $\cI_{1,2}$ that join in the first quadrant making a cusp;
the point of their intersection is the nearest (w.r.t.\ the Euclidean
metric on~$0z_1z_2$) inflection point located
on the spirals that unroll from the origin.

The first component~$\cI_1$ consists of the inflection points where
the trajectories start to repulse from the separatrix~$y_2(x)$
(the unstable solution $y_2$ tends from above to the ray~$y=x$, see
Proposition~\ref{AsymptSepState} below) and turn towards their local
maxima on the ray~$y=x/2$. At the infinity $x^2+y^2\to\infty$, the
component~$\cI_1$ approaches the diagonal~$y=x$ from below.
   % Если вообще -- продолжается бесконечно!
%
The second component~$\cI_2$ goes left from its intersection with
$\cI_1$ and contains the inflection point of the separatrix. Then the
curve~$\cI_2$ enters the second quadrant. There, it describes the
moments when the phase curves, having crossed the $0x$ axis
at~$x<0$, repulse from the separatrix and turn upward, possessing the
vertical asymptotes at infinity.
\end{rem}

\subsection{}
Having achieved the inflection point, the trajectories of
system~\eqref{xysyst} approach the infinity following
one of the three schemes below. The expansions for the solutions
$\psi_1(\tau)$, $\psi_1^*(r)$, and
$\psi_2(\tau)$ of the Riccati equation~\eqref{Riccati},
which were obtained in the previous section, determine the asymptotes
for the phase curves of all the three types.

\begin{state}   %%%[The asymptote of the separatrix~$y_2$]
   \label{AsymptSepState}  %=7
The diagonal ${y=x}$ is the asy\-m\-p\-tota of the separatrix~$y_2(x)$ in a
neighbourhood of infinity\textup{;} the curve $y=y_2(x)$ approaches the
diagonal from above.
\end{state}

\begin{proof}
Using expansion~\eqref{psi2}
%$$ \tag{ } %
%\psi_2(\tau)=1-\sqrt2+\frac{2\sqrt2(\sqrt2-1)}{3\tau}+O(\frac1{\tau^2})
%   \tag{\ref{psi2}}
%$$ %\end{equation}
for the solution $\psi_2(\tau)$ as $\tau=4r^2\to\infty$,
and taking into account the transformation from Eq.~\eqref{xysyst}
to Eq.~\eqref{z_syst}, we obtain
\begin{multline*}
y-x=2z_2=2r\sin\varphi
=2r\sin\bigl((\varphi-\tfrac{\pi}8)+\tfrac{\pi}8\bigr)=\\
=2r\cos\tfrac{\pi}8\cos(\varphi-
  \tfrac{\pi}8)\cdot
  \bigl(\tan(\varphi-\tfrac{\pi}8)+\tan\tfrac{\pi}8\bigr)=
2r\cos\tfrac{\pi}8
 \bigl(\psi_2(\tau)+\sqrt2-1\bigr)\bigr/\sqrt{1+\psi_2^2(\tau)}=\\
=2r\cos\tfrac{\pi}8\Bigl(\frac{2\sqrt2(\sqrt2-1)}{12r^2}
 +O\bigl(\frac{1}{r^4}\bigr)
  \Bigr) \bigl/ \sqrt{1+\psi_2^2(\tau)}
  \xrightarrow{r\to\infty} +0.
\end{multline*}
This argument concludes the proof.
\end{proof}

The separatrix divides the trajectories in the first (consequently, in
the third) quadrant in two classes. The curves~$y_1^*$ of the first
type repulse from the separatrix
and go vertically upward, not leaving the
quadrant. Other trajectories dive under the separatrix and achieve the
extrema on the ray~$y=x/2$. Then they  intersect the axis~$0x$, whence
they either turn down and have a vertical asymptote at~$-\infty$
or reach the third quadrant.
Then the whole situation is reproduced up to the central symmetry.

Now we give a rigorous proof of these properties.

\begin{state}   %%%[The asymptote of the solution~$y^*_1$]
   \label{Y1StarAsympt}   %=8
The coordinate axis~$x=0$ is the vertical asymptote for the
solution~$y_1^*$, which approaches it from the left.
\end{state}

\begin{proof}
Using expansion~\eqref{psi1star}
%$$
%\psi_1(\tau)=\sqrt2-1+\frac{2\sqrt2(\sqrt2-1)}{\tau}+O(\frac1{\tau^2}).
%$$
and inverting the transformation $0xy\mapsto 0z_1z_2$,
we let $r\to\infty$ and hence obtain
\begin{multline}
x=2(z_1-z_2)=2r(\cos\varphi-\sin\varphi)=\\
=-2\sqrt2r\sin(\varphi-\tfrac{\pi}4)=
-2\sqrt2r\sin\bigl((\varphi-\tfrac{\pi}8)-\tfrac{\pi}8\bigr)=\\
=-2\sqrt2\cos\tfrac{\pi}8\cdot
r\bigl(\tan(\varphi-\tfrac{\pi}8)-\tan\tfrac{\pi}8\bigr)
  \cos(\varphi-\tfrac{\pi}8)=\\
=-2\sqrt2\cos\tfrac{\pi}8\cdot
r\bigl(\psi_1(\tau)-\sqrt2+1\bigr)\bigr/\sqrt{1+\psi_1^2(\tau)}=\\
=-2\sqrt2\cos\tfrac{\pi}8\cdot
r\Bigl(
\frac{2\sqrt2(\sqrt2-1)}{4r^2}+O\bigl(\frac1{r^4}\bigr)
\Bigr)\bigl/\sqrt{1+\psi_1^2(\tau)} \to-0.
\label{xTendsToZero}   %=25
\end{multline}
According to~\eqref{radial}, the function~$r$ grows infinitely along the
trajectories. Consequently, the function $R(t)=\sqrt{x^2(t)+y^2(t)}$
also tends to $+\infty$ along the solutions $\{x(t)$,\ $y(t)\}$.
By~\eqref{xTendsToZero}, the coordinate~$x$ tends to
zero from the left, therefore the axis~$x=0$ is the vertical asymptote
of the trajectory~$y_1^*(x)$.
\end{proof}

%%%%%%%%%%%%%%%%%%%%%%%%%%%%%%%%%%%%%%%%%%%%%%%%%%%%%%%%%%%%%%%%%%%%%%%%
Let us analyze the asymptotic behaviour of the trajectories~$y_1(x)$,
which correspond to the solutions $\psi_1(r)$ of the Riccati
equation~\eqref{Riccati}.
Using expansion~\eqref{psi1} for an arbitrary solution~$\psi_1(r)$,
we describe the asymptotes on the plane~$0xy$.
We recall that the choice~$C>0$ in formula~\eqref{psi1} corresponds to
solutions located above the graph $\psi_1^*(\tau)$.
The expansion of the function $\psi_1^*$ itself is given at $C=0$,
and the asymptote of its image on the plane $0xy$ was obtained in
Proposition~\ref{Y1StarAsympt}.
If $C<0$, then the solutions $\psi_1(r)$ are located
between~$\psi_1^*(\tau)$ and~$\psi_2(\tau)$.

Consider the phase plane~$0xy$. The solutions of the class $y_1$ grow
infinitely in the second quadrant if~$C>0$. Suppose~$C<0$,
then the representatives of this class go to infinity
between the axis~$0y$ and the diagonal~$y=x$ in the first quadrant;
recall that the diagonal is the asymptote of the separatrix~$y_2$.
Using an appropriate modification of the
proof of Proposition~\ref{Y1StarAsympt}, we obtain the estimate
\begin{multline*}
x=-\frac{2\sqrt2r\cos\frac{\pi}8}{\sqrt{1+u^2(r,C)}}
\bigl(u(r,C)-\tan\tfrac{\pi}8\bigr)   %=\\
=-\frac{2\sqrt2r\cos\frac{\pi}8}{\sqrt{1+u^2(r,C)}}
\left(\frac{C}{r}+\frac{D}{r^2}+O\bigl(\frac1{r^3}\bigr)\right),
\end{multline*}
where
  %$$ \tag{ } %\label{PutD}
$D=\tfrac{1}{2}\,[ {2-\sqrt2+(\sqrt2+1)\cdot C^2} ]$
  %$$ %\end{equation}
and expansion~\eqref{psi1} is substituted for~$u(r$,\ $C)$.
Hence we finally obtain

\begin{state}\label{GeneralY1State}
All solutions of the class~$y_1$ have vertical asymptotes,
$$
\lim_{r\to\infty}x(r,C)=-2\sqrt2C\cos^2\tfrac{\pi}8=-(\sqrt2+1)C.
$$
The graphs of solutions contained in the upper half\/-\/plane
approach these asymptotes from the left.
\end{state}

\begin{rem}
Without loss of generality, let~$C<0$. Then the vertical straight line
$x=-(\sqrt{2}+1)\cdot C$, which is located in the right half\/-\/plane,
is the vertical asymptote for two trajectories of system~\eqref{xysyst}.
If a solution~$y_1$ is greater than the separatrix~$y_2$, then it
tends to the asymptote from the \emph{left}.
Another trajectory dives under the separatrix and approaches
the same asymptote from the \emph{right} in the fourth quadrant.
The constants $C>0$ describe centrally
symmetric curves in the third and second quadrants, respectively.
\end{rem}

\begin{rem}\label{CompareWithLimit}
From Propositions~\ref{AsymptSepState}--\ref{GeneralY1State}
it follows that for large~$R$ solutions of
Eq.~\eqref{phase_crv} are approximated by the exact solutions
of limit equation~\eqref{phase_crv3}. The constants~$C$ in~\eqref{psi1}
and~$\delta$ in~\eqref{ApproxDelta} are related by the formula
$(\sqrt{2}+1)\cdot C=\pm 1/\sqrt{2\delta}$.
\end{rem}

%%%%%%%%%%%%%%%%%%%%%%%%%%%%%%%%%%%%%%%%%%%%%%%%%%%%%%%%%%%%%%%%%%
\section{The spiral minimal surfaces}\label{SecMinSurf}
By construction, solutions~$h(q)$ of Eq.~\eqref{BilaCubic}
determine the profiles $\Pi=\{z=h(q)\mid\rho=1$,\ $q\in\BBR\}$
of the minimal surfaces~$\varSigma\subset\BBE^3$
on the cylinder $Q=\{\rho=1\}$; we thus have~$\Pi=\varSigma\cap Q$.
We further recall that a minimal surface in~$\BBE^3$ is extended from
the profile~$\Pi$ in agreement with Corollary~\ref{DrawBilaSpiral}.
Yet, not each surface in~$\BBE^3$
assigned to a solution of~\eqref{BilaCubic} is \emph{minimal}.
Let us study this aspect in more detail.

\subsection{The selection rule}
In this subsection we formulate a rule that defines which components of
the phase trajectories for Eq.~\eqref{BilaCubic} provide minimal
surfaces.    %%%~\eqref{xysyst},
Using this rule, we conclude that the graphs $\{z=\chi(\xi$,\ $\eta)\}$
of solutions $\chi=\rho\cdot h(\theta-\ln\rho)$ attach nontrivially to
each other.

Let us recall a well\/-\/known property of the minimal
surfaces~\cite{Nitsche}. A smooth
two\/-\/dimensional surface in space is
minimal iff at any point its mean curvature~$H$ vanishes.
Hence each point of the surface is a saddle
(of course, we assume that the surface at hand is not a plane).

Further on, we denote by~$\frac{\dd}{\dd\theta}$ the respective
coordinate vector attached to a point of the cylinder~$Q$.
   %%%Looking onto the plane~$z=0$ from above ($z\to+\infty$),
We see that the logarithmic spirals $q=\theta-\ln\rho=\const$
are convex, and at any point of the surface
the nonzero curvature vector
has a positive projection onto $\frac{\dd}{\dd\theta}$.
Therefore we request that the curvature vectors at all points
of the profiles~$\Pi=\{z=h(q)$, $\rho=1\}\subset Q$
must have negative projections on~$\frac{\dd}{\dd\theta}$.

The above condition upon solutions~$h(q)$ of Eq.~\eqref{BilaCubic} can
be reformulated as follows: $h'\cdot h''>0$.
This inequality is the rule for selecting the components of the phase
curves on the plane~$0xy$. Namely, consider a
neighbourhood~$\Omega_\Gamma$ of a point~$(x,y)$ on a phase
curve~$\Gamma$. The minimal surface~$\varSigma_\Gamma$ is
assigned to this component of the curve if the coordinates
$x$ and $y$ satisfy the inequality
$$
y^3\cdot(y^3-(x-y)^3+2y-x)>0.
$$
Hence we obtain~$y>x/2$ for~$y>0$ and~$y<x/2$ whenever~$y<0$, see the
proof of Proposition~\ref{ExtremaState}.
If not all points of a phase curve satisfy this inequality
(actually, each curve while it stays near the focus has to pass through
the domains where it is not satisfied), then several components of the
profiles~$\Pi$ and several minimal surfaces~$\varSigma$
are assigned to this curve.
A continuous motion along the phase trajectory~$\Gamma$
corresponds to different (possibly, distant from each other)
components of the profile~$\Pi$.

In Sec.~\ref{SecPhasePlane} we described two classes of
trajectories on the plane~$0xy$, the curves~$y_1(x)$ with vertical
asymptotes and the separatrix~$y_2(x)$ with the diagonal asymptote~$y=x$.
Clearly, a small neighbourhood of the focus at the origin
corresponds to oscillations of any profile, and some parts of the
oscillating curves $z=h(q)$ are prohibited by the selection rule.
At infinity, the solutions of first type
describe the finite size fragments of the
profiles~$\Pi_1=\{z=h(q)\}$; these fragments have the vertical tangent
at a finite value of the function~$h$.
Conversely, when the separatrix approaches the diagonal asymptote, it
defines the exponentially growing profile~$\Pi_2$ that turns
infinitely many times around the cylinder~$Q$.

\begin{example}[\cite{YS2005}]\label{ExYS2005}
Let us choose the constant~$C<0$ such that $y_x^*(-(\sqrt{2}+1)C)=0$.
Consider the following
two components of the phase curves~$\Gamma=\{y=y_1(x$,\ $C)\}$:
\begin{align*}
\Gamma_{-\infty}&=\{(x,y)\mid -(\sqrt{2}+1)C \geq x(t)>0,\
     0\geq y=y_1^*(x)\xrightarrow{x\to+0}-\infty\},\\
\Gamma_{+\infty}&=\{(x,y)\mid 0\leq x(t)< -(\sqrt{2}+1)C,\ %
     y=y_1(x)\xrightarrow[x\to-(\sqrt{2}+1)C]{}+\infty\}.
\end{align*}
Also, let us use the focus~$(0,0)$, which corresponds to the trivial
solution~$h\equiv0$ and hence amounts to the plane~$z=0$ in~${\BBE}^3$.
Using these curves, we construct the closed
profiles~$\Pi_1(C)\subset Q$, see Fig.~1.   %%% \cite{YS2005}
\[
{
\unitlength=1mm
\begin{picture}(50,25)(0,3.5)
\put(5,5){\llap{Fig.\,1}}
\put(0,15){\vector(1,0){50}}
\put(10,0){\vector(0,1){30}}
\put(49,12){$q$}
\put(9,27){\llap{$h$}}
\put(10,15){\circle*{1}}
\put(34,15){\circle*{1}}
\put(9,16){\llap{$R$}}
\put(35,16){$S$}
\put(18.5,21){\llap{$\sigma_1$}}
\put(18.5,7){\llap{$-\sigma_1$}}
\put(35,21){$\sigma_1^*$}
\put(35,7){$-\sigma_1^*$}
\qbezier(10,15)(22,19)(22,27)
   %\put(10,27){\oval(24,24)[br]}
\qbezier(10,15)(22,11)(22,3)
   %\put(10,3){\oval(24,24)[tr]}
\put(22,15){\oval(24,24)[r]}
\put(22,16){$\Pi$}
\end{picture}
\phantom{MMMMMMMM}
\begin{picture}(50,25)(0,3.5)
\put(-2,5){\llap{Fig.\,2}}
\put(0,2){\vector(1,0){42}}
\put(2,0){\vector(0,1){30}}
\put(42.5,3){$q$}
\put(1,27){\llap{$h$}}
\put(14,10){\llap{$\sigma_1$}}
\put(24,20){\llap{$\tilde{\sigma}_1$}}
\put(28,14){$\sigma_2$}
\qbezier(10,6)(28,9)(30,26)
   %\put(6,26){\oval(40,40)[br]}
\qbezier(10,6)(20,10)(20,16)
\qbezier(20,16)(28,18)(30,26)
\end{picture}
}
\]
By Corollary~\ref{DrawBilaSpiral}, the points $R$ and $S$ correspond to
the logarighmic spirals on the plane~$z=0$ in~$\BBE^3$.
Shifting the point~$R$ and the components~$\pm\sigma_1$ towards~$S$,
we make the length $|RS|$ of the profile~$\Pi$ comparable with~$2\pi$.
Attaching the respective number of the profiles one after another
on the cylinder~$Q$, we place the edge~$R$ of each profile on
the spiral~$S$ of the previous profile. Hence the resulting minimal surface
becomes self\/-\/supporting.
\end{example}

\begin{example}
Consider the following parts of the phase curves:
a solution~$y_1(x)$ as it tends to~$+\infty$,
another solution~$\tilde{y}_1(x)$ of this class located
\emph{below} the separatrix before the extremum on~$y=x/2$,
and the separatrix~$y_2(x)$ itself.
Then we obtain the closed profile with finite cross\/-\/section
in the upper half\/-\/plane~$h>0$.
Indeed, we attach the corresponding curves $\sigma_1$,
$\tilde{\sigma}_1$, and $\sigma_2$ as shown on Fig.~2.
\end{example}

\subsection{The profiles $\Pi=\{z=h_{1,2}(q)\}$ and their approximations}
 %In this section we consider solutions~$h(q)$ of Eq.~\eqref{BilaCubic}.
 %By construction, these solutions determine the profiles~$\Pi=\{z=h(q)\}$
 %of the minimal surfaces~$\varSigma$ on the
 %cylinder~$Q\subset\BBE^3$, see Corollary~\ref{DrawBilaSpiral} on
 %p.~\pageref{DrawBilaSpiral}.
Now we obtain the asymptotic estimates for the profiles $\Pi_1$.
In this subsection we assume that the absolute value of the derivative
$h'\equiv\tfrac{dh}{dq}$ (and, possibly, of the function~$h$ itself)
is large.

Let the constant~$C$ in~\eqref{psi1} be arbitrary and put
\[
D=\tfrac{1}{2}\bigl[2-\sqrt{2}+C^2\cdot(\sqrt{2}+1)\bigr]
\]
as before. Again, we use expansion~\eqref{psi1} as~$r\to\infty$.
Then we get the equivalence~($\sim$)
\begin{multline*}
h=x=2(z_1-z_2)=2r(\cos\varphi-\sin\varphi)=\\
=\frac{2\sqrt2\cos\frac{\pi}8\cdot
r(\sqrt2-1-u(r,C))}{\sqrt{1+u^2(r,C)}}
\sim -(\sqrt2+1)\left(C+\frac{D}{r}\right),
\end{multline*}
which holds up to~$O(\tfrac{1}{r^2})$.
The derivative~$h'$ is expressed through a solution of
Eq.~\eqref{Riccati} in the following way:
\begin{multline*}
h'=y=2r\cos\varphi=2r\cos\tfrac{\pi}8\cos\bigl(\varphi-\tfrac{\pi}8\bigr)
\left(1-\tan\tfrac{\pi}8\tan\bigl(\varphi-\tfrac{\pi}8\bigr)\right)=\\
=\frac{2r\cos\frac{\pi}8(1-(\sqrt2-1)u(r,C))}{\sqrt{1+u^2(r,C)}}.
\end{multline*}
As~$r\to\infty$, the above formula is equivalent ($\sim$) to
\begin{multline*}
2r\cos^2\tfrac{\pi}8\left(1-(\sqrt2-1)^2-(\sqrt2-1)
\Bigl(\frac{C}{r}+\frac{D}{r^2}\Bigr)\right)=\\
=r\Bigl(1+\frac{\sqrt2}2\Bigr)(2\sqrt2-2)-
\Bigl(1+\frac{\sqrt2}2\Bigr)(\sqrt2-1)\Bigl(C+\frac{D}{r}\Bigr)=\\
=\sqrt2r-\frac{\sqrt2}2\Bigl(C+\frac{D}{r}\Bigr).
\end{multline*}
This reasoning shows that for large~$r$ we have
$$
h\sim-(\sqrt2+1)\left(C+\frac{D}{r}\right),\quad
h'\sim\sqrt2\,r-\frac{\sqrt2}2\left(C+\frac{D}{r}\right).
$$
Eliminating~$r$ from these formulas, we obtain the differential
equation
$$  %$$ %\label{four}
h'=\frac{\sqrt2}2\frac{h}{\sqrt2+1}
  -\frac{\sqrt2D}{C+\frac{\displaystyle h}{\displaystyle 1+\sqrt2}}.
$$  %$$ %\end{equation}
By definition, put~$g=\frac{h}{\sqrt2+1}$. Then we get
\begin{equation}\label{five}   %=28
(\sqrt2+1)g'=\frac{\sqrt2}2g-\frac{\sqrt2D}{C+g}.
\end{equation}
Let us consider two cases: $C=0$ and $C\ne0$.
If~$C=0$, then we see that
$$
\frac{(\sqrt2+1)2g\,dg}{\sqrt2(g^2-(2-\sqrt2))}=dq
\quad\Longrightarrow\quad
\frac{\sqrt2+1}{\sqrt2}\ln|g^2-(2-\sqrt2)|=q+A,
$$
where $A=\const$. Resolving this equality w.r.t.~$g$, we obtain
$$
g=\pm\sqrt{2-\sqrt2+\exp(\frac{\sqrt{2}}{\sqrt2+1}q+A)}.
$$
Now let~$C\ne0$; then we express~$g(q)$ from Eq.~\eqref{five}.
Therefore we get ${h=(\sqrt{2}+1)\,g(q)}$ in the form of an integral.
   %%%Hence we obtain the assertion.

\begin{state}\label{AsymptH1Star} %=10
Suppose~$|h'|\to\infty$.
If~$C=0$, then the asymptotic behaviour of the solution~$h_1^*(q)$ of
Eq.~\eqref{BilaCubic}, which is assigned to the phase curve~$y_1^*(x)$,
is described by the formula
$$
h_1^*(q)\sim
\pm\sqrt{2+\sqrt2+(\sqrt2+1)^2\exp(\sqrt2(\sqrt2-1)q+A)},
$$
here~$A$ is an arbitrary constant.
If~$C\ne0$, then the approximation for~$h_1^*(q)$ is given implicitly
through the integral for Eq.~\eqref{five},
\begin{multline*}
{\bigl| (\sqrt{2}-1)h+\tfrac{1}{2}C-B\bigr|}^{2B+C}\cdot
{\bigl| (\sqrt{2}-1)h+\tfrac{1}{2}C+B\bigr|}^{2B-C} %={}\\
{}=A\cdot\exp\bigl(2(2-\sqrt{2})\,Bq\bigr),
\end{multline*}
here $A>0$,  $B=\sqrt{\tfrac{1}{4}C^2+2D}$,
$C\in\BBR$, and
$D=\tfrac{1}{2}\bigl[2-\sqrt{2}+C^2\cdot(\sqrt{2}+1)\bigr]$.
\end{state}

Finally, we use Proposition~\ref{AsymptSepState} and describe
the profile~$\Pi_2$ assigned to the separatrix~$y_2$.
The asymptote~$y=x$ defines the equivalence~$h'\sim h$
as~$|h|$,\ $|h'|\to\infty$. Consequently, the profile~$\Pi_2$
is approximated by solutions of the differential equation~$h'=h$.

\begin{state}\label{ExpForSep}   %=11
Let $h$ be large, then the function~$h_2(q)$ grows exponentially
in~$q\in\BBR$\textup{:}
$h_2\sim A\cdot\exp(q)$, where~$A=\const$.
\end{state}

\subsection{Approximations of the minimal surfaces}
From Propositions~\ref{AsymptH1Star} and~\ref{ExpForSep} we obtain the
following assertion.

\begin{theor}\label{AsymptTh}
Let $\rho$,\ $\theta$ be the polar coordinates on the plane~$0\xi\eta$
and put~$q=\theta-\ln\rho$, $h={z\bigr|}_{\rho=1}$.
\begin{enumerate}
\item  %{1.}
For large values of $h'(q)$, there is a $\phi$-\/invariant minimal
surface $\varSigma_1^*=\{z=\rho\cdot h_1^*(\theta-\ln\rho)\}$,
which is approximated by the graph of the function
\begin{equation}\label{AsymptSigma1Star}  %=29
z=\rho\cdot\sqrt{2+\sqrt2+
\frac{(\sqrt2+1)^2\exp(\sqrt2(\sqrt2-1)\theta+A)}{\rho^{\sqrt2(\sqrt2-1)}}},
\quad A=\const.
\end{equation}
Also, there is a class of the $\phi$-\/invariant minimal surfaces
$\varSigma_1$ that correspond to the asymptotic
formulas for~$h_1(q)$ in Proposition~\textup{\ref{AsymptH1Star}}.
This class contains $\varSigma_1^*$ as a particular case.
\item   %{2.}
Suppose that the absolute values~$|h_2|$ are sufficiently large such
that the phase curve~$h_2'(h_2)$ is near to the diagonal~$h'=h$.
Then there is the $\phi$-\/invariant minimal surface
$\varSigma_2=\{z=\rho\cdot h_2(\theta-\ln\rho)\}$
whose approximation is given by the graph of the function
\begin{equation}\label{AsymptSigma2}   %=30
z = A\,\rho\cdot\exp(\theta-\ln\rho)=A\cdot\exp(\theta),\quad A=\const.
\end{equation}
\end{enumerate}
The constant~$A\in\BBR$ determines the rotation of the
surfaces around the axis~$0z$. The choice~$A<0$ in~\eqref{AsymptSigma2}
defines the reflection symmetry~$z\mapsto-z$.
Both surfaces $\varSigma_1^*$, $\varSigma_{2}\subset{\BBE}^3$
have a singular point at the origin.
\end{theor}

Finally, let us plot the graph of a spiral minimal surface
determined by the profiles~$\Pi$ on Fig.~1, see Example~\ref{ExYS2005}.
\[
{
\unitlength=1mm
\begin{picture}(70,38)(0,3.5)
\put(0,5){\llap{Fig.\,3}}
   %\put(10,0){\vector(0,1){30}}
   %%% No.1
\qbezier(0,15)(12,19)(12,27)
\qbezier(0,15)(12,11)(12,3)
\qbezier(12,27)(24,27)(24,15)
\qbezier(24,15)(24,3)(12,3)
   %\put(12,15){\oval(24,24)[r]}
  \qbezier(12,27)(12,43)(36,43)
  \qbezier(36,43)(60,43)(60,27)
      %\put(36,27){\oval(48,32)[t]}
  \qbezier(0,15)(0,21)(12,27)
  \qbezier(24,15)(24,27)(36,27)
   %%% No.2
\qbezier(60,27)(60,19)(48,19)
\qbezier(48,19)(36,19)(36,27)
      %\put(48,27){\oval(24,16)[b]}
\qbezier(60,14)(60,9)(48,9)
\qbezier(48,9)(36,9)(36,14)
      %\put(48,15){\oval(24,16)[b]}
\qbezier(66,21)(66,13)(48,13)
\qbezier(48,13)(30,13)(30,21)
      %\put(48,21){\oval(36,24)[b]}
\qbezier(66,21)(60,25)(60,27)
\qbezier(30,21)(36,25)(36,27)
   %\put(10,27){\oval(24,24)[br]}
   %\put(10,3){\oval(24,24)[tr]}
   \qbezier(66,21)(66,24)(60,27)
   \qbezier(30,21)(30,24)(36,27)
     \qbezier(48,27)(54,27)(54,20)
     \qbezier(60,27)(54,27)(54,20)
       %\qbezier(54,21)(54,19)(52,17)
     \qbezier(36,27)(40,26)(40,21)
     \qbezier(40,21)(40,25.5)(40.5,25.5)
   %%% No.3
\qbezier(36,27)(36,31)(42,31)
\qbezier(42,31)(48,31)(48,27)
   \qbezier(49.5,25.5)(48,26)(48,27)
     %\qbezier(49.5,25.5)(49.5,26)(48,27)
   \qbezier(49.5,25.5)(49.5,23.5)(45,23.5)
   \qbezier(45,23.5)(40.5,23.5)(40.5,25.5)
       %\qbezier(40.5,25.5)(40.5,27)(42,27)
     \qbezier(40.5,25.5)(42,25.5)(42,27)
     \qbezier(45,27)(46.5,27)(46.5,25.2)
     \qbezier(48,27)(46.5,27)(46.5,25.2)
   %%% No.4
\qbezier(48,27)(48,25)(45,25)
\qbezier(45,25)(42,25)(42,27)
   \qbezier(48,24)(48,22.5)(45,22.5)
   \qbezier(45,22.5)(42,22.5)(42,24)
   %%% No.5
\qbezier(42,27)(42,28)(43.5,28)
\qbezier(43.5,28)(45,28)(45,27)
   %%% No.6
\put(44,27){\oval(2,2)[b]}
\put(44,27){\oval(3.5,3)[b]}
   %%%\put(45,25.7){\oval(2.6,2.6)[tl]}
   %%%
\put(43.5,27){\circle*{1.2}}
\end{picture}
}
\]

We conjecture that the spiral minimal surfaces may be
relevant in astrophysics owing to their visual similarity with the
spiral galaxy\/-\/like objects. Also, we indicate their application in
fluid dynamics: these surfaces can realize the minimum of the free
energy~$\mathcal{F}$ of a vortex.
   %%%Finally, these surfaces may prove useful in architecture.

%%%%%%%%%%%%%%%%%%%%%%%%%%%%%%%%%%%%%%%%%%%%%%%%%%%%%%%%%%%%%%%%%%%
\subsection*{Acknowledgements}
The authors thank R.\,Vitolo and D.\,V.\,Pelinovski\v{\i}
for helpful discussions.
A part of this research was carried out while A.\,K. was
visiting at the University of Lecce.
A.\,K. was partially supported by University of Lecce
grant no.650~CP/D.

\end{document}